\documentclass[11pt]{amsart}
\usepackage{amsfonts}
\usepackage{amsmath,amsthm}
\usepackage{bbm}

\usepackage{latexsym}
\usepackage{latexsym}
\usepackage{array}
\usepackage{amssymb}
\usepackage{enumerate}
\usepackage{comment}
\usepackage{graphicx}
\usepackage{tikz}
\usepackage{color}
\usepackage{bbm}

\DeclareFontFamily{U}{wncy}{}
\DeclareFontShape{U}{wncy}{m}{n}{<->wncyr10}{}
\DeclareSymbolFont{mcy}{U}{wncy}{m}{n}

\usepackage[english]{babel}
\usepackage{fdsymbol}
\usepackage{color}
\usepackage[latin1]{inputenc}

\numberwithin{equation}{section}

\theoremstyle{plain}
\newtheorem{theorem}{Theorem}[section]
\newtheorem*{theorem*}{Theorem}
\newtheorem{lemma}[theorem]{Lemma}
\newtheorem{proposition}[theorem]{Proposition}

\theoremstyle{remark}

\newtheorem*{lem*}{Lemma}
\newtheorem*{sublem*}{Sublemma}
\newtheorem*{remark*}{Remark}
\newtheorem*{NB*}{NB}

\newcommand{\supp}{\mathop{\rm supp}\nolimits}

\newcommand{\R}{ \mathbb{R} }

\newcommand{\C}{ \mathbb{C} }
\newcommand{\Z}{ \mathbb{Z} }

\newcommand{\K}{ \mathbb{K} }

\newcommand{\T}{ \mathbb{T} }
\newcommand{\A}{ \mathbb{A} }

\newcommand{\vs}{ \vec s }

\newcommand{\cA}{ \mathcal{A} }

\newcommand{\cF}{ \mathcal{F} }

\newcommand{\EE}{ {\mathbb E}}

\newcommand{\cZ}{ \mathcal{Z} }

\newcommand{\om}{ \omega }

\newcommand{\ga}{\gamma }

\renewcommand{\phi}{ \varphi }
\newcommand{\oms}{ \omega^{12}_{3s}  }

\newcommand{\eps}{\varepsilon}

\newcommand{\al}{ \alpha }
\newcommand{\zz}{\mathfrak z}

\newcommand{\la}{ \lambda }

\newcommand{\dess}{\delta'^{12}_{3s}}

\newcommand{\tauz}{\tau_0}

\newcommand{\be}{\begin{equation}}
\newcommand{\ee}{\end{equation}}
\newcommand{\ben}{\begin{equation*}}
\newcommand{\een}{\end{equation*}}

\newcommand{\lan}{ \langle }

\newcommand{\ran}{ \rangle}

\newcommand{\p}{ \partial}

\newcommand{\e}{ \text{e} }
\newcommand{\ai}{ \mathfrak{a} }

\newcommand{\lbl}{\label}

\newcommand{\qu}{\quad}
\newcommand{\qmb}{\quad\mbox}
\newcommand{\qnd}{\qmb{and}\qu}

\newcommand{\no}{\mathfrak{n}}

\newcommand{\dep}{\delta'^{12}_{3s}}

\title[Wave turbulence]
{An algebraic geometry question and  wave turbulence}

\author{Sergei  Vl\u adu\c t}
\address{Sergei  Vl\u adu\c t \\ Aix Marseille Universit\'e, CNRS,  I2M UMR 7373, 13453, Marseille,
France } \email{ vladut@iml.univ-mrs.fr, sergevladuts@ya.ru}

\begin{document}

\begin{abstract} This paper is a sequel to \cite{seclim}.
We formulate a natural algebraic geometry conjecture,
give some of its number theoretic and analytical consequences,
and show that those can be used to get further advances in wave turbulence theory.
\end{abstract}

\date{}

\maketitle
\section{Introduction}\label{s1}
This note (a sequel to \cite{seclim}) is devoted to stating an elementary formulated conjecture on the dimension of
the intersection of some quadrics in higher dimensional (affine) space. Following \cite{seclim} we deduce from it some interesting arithmetical and analytical consequences and show that those consequences can be used to get some new advances in the Zakharov--L'vov theory of  wave turbulence. This note does not contain essentially new results, but stresses that one can  get rather important advances in  Zakharov--L'vov  theory by calculating the dimension  of some simply looking  affine algebraic sets. 

I address this note especially to algebraic geometers in a hope  that somebody  will be able to calculate this dimension confirming  our main conjecture and, therefore, its numerous prediction of arithmetic, analytic and PDE nature.

Our paper is organized as follows: 

-- in Section 2 we formulate the conjecture and give some (very partial) results in its direction; in particular, we prove two simplest cases of the conjecture, already used in \cite{seclim}.

-- Section 3 contains some arithmetic and analytic consequences of the conjecture; in Sections 2 and 3 our  presentation is rather complete and self-contained.

-- On the contrary, Sections 4 and 5 give a very brief and incomplete exposition of applications to the Zakharov--L'vov theory. We address the reader to \cite{seclim} and  references therein for more detailed information. However, I hope that the reader will have an idea how to  apply the results  of Section 3 to wave turbulence problems.

-- Section 4 contains a brief description of the  Zakharov--L'vov  model(s),  which essentially  reduce(s) difficult physical problems to (still difficult) questions on the asymptotic behavior of solutions to  non-linear stochastic PDE's.

-- In Section 5 we describe an application of the results of Sections 2 and 3 to the corresponding asymptotic problems.

\medskip {\em Acknowledgment.} I thank A.Maiocchi for  sharing with me his note \cite{Mai}. 
\section{The question}
 Let $\K$ be an algebraically closed field, let $N\ge 2,d\ge 2$ be two fixed integers, and let $\al=(\al_{ij})$, be an $N\times N$ skew-symmetric matrix whose elements belong to the set $\{-1, 0, 1\}\subset \K$, 
 without zero lines and rows. 
Consider a family of quadratic forms on $\K^{dN}$
\be\label{om}
\om_j(z)=z_j\cdot \sum_{i=1}^N \al_{ji}z_i, \qquad 1\leq j\leq N,
\ee 
where  
$
z=(z_1,\ldots,z_N),
$
$z_j\in\K^d,$ and $x\cdot y$ denotes the standard scalar product $ \sum_{i=1}^{dN} x_iy_i$.
Let us  define the
 geometric quadrics $Q_j:=\{z\in \K^{Nd}:\, \omega_j(z)=0\}$ and consider their intersection
 $$Q=\cap_{j=1}^N Q_j\subset \K^{dN}=\A_\K^{dN} =\A^{dN} $$ which is an affine algebraic set (the affine cone over a projective algebraic set). 
  Note that $Q=\bigcap_{j=1}^{N-1} Q_j$ since the skew symmetry of the matrix $\al$ implies 
  that $\om_1+\ldots+\om_N=0$.
 
By suitably rearranging  indices $i$  and possibly
multiplying $\om_i$ by $-1$, $\omega_1$ may be written in the form 
$
\omega_1(z) =z_1\cdot \sum_i \alpha_{1i} z_i
$ with $\alpha_{1N}=1$ which we assume below. Define $v= \sum_i\alpha_{1i} z_i$ so that 
\be\lbl{z->v}
\omega_1(z) =z_1\cdot v \qnd z_N =
v-\sum_{1<i<N}\alpha_{1i}z_i\,.
\ee
If we fix  $(z_1,v)\in \K^{2d}$  the remaining quadratic forms
$\omega_j$ with $1<j<N$ as functions of $(z_2,\ldots,z_{N-1})\in \K^{(N-2)d}$ 
 become polynomials of degree at most two, with no
constant term. Namely
\begin{equation}
 \label{q_tilde}
q_j(z_2,\ldots,z_{N-1};z_1,v)=z_j\cdot\Big(\alpha_{j1} z_1 +
\alpha_{jN} v+
\sum_{1<i<N}(\alpha_{ji}-\alpha_{jN}\alpha_{1i})z_i\Big)\,.
\end{equation}

Consider the sets 
$$
\tilde Q_j(z_1,v)=\{(z_2, \dots, z_{N-1}) :q_j=0\} \subset \A^{(N-2)d}, \quad 1<j< N, 
$$
and their intersection $ Q(z_1,v)=\cap_{1<j<N}\tilde Q_j(z_1,v)$. 
Slightly abusing notation, below we write $
Q_1=\{(z_1,v)\in\A^{2d}:z_1\cdot v=0\}.$  

Therefore, we have the linear projection
$$ \pi=\pi_v:\A^{Nd}\longrightarrow \A^{2d},\; \pi (z_1,\ldots,z_N)=(z_1,v= \sum_{1<i\le N}\alpha_{1i}z_i),$$
for which
$$\pi(Q)=Q_1,\;\pi^{-1}(z_1,v)=Q(z_1,v). $$

 Let $$A_2:=\{(z_1,v)\in \A^{2d}:z_1\ne 0,v\ne 0. \} $$
 
 Let us formulate our main elementary algebraic geometry question, which we pose as a conjecture.\bigskip
 
 {\bf Main conjecture.} (MC) {\em For any  $(z_1,v)\in Q_1\cap A_2$ we have}
 \be\label{mc} \dim Q(z_1,v)=(N-2)d.\ee

 In other words, the quadrics $\{\tilde Q_j(z_1,v),1<j<N\}$ intersect transversally.\bigskip

  Note that if $\al=\al_1\oplus\al_2$ and the conjecture holds for both $\al_1$ and $\al_2,$ it holds for $\al$
 and thus discussing  MC we can and will suppose that $\al$ is irreducible; also the condition  $(z_1,v)\in Q_1\cap A_2$ is supposed to hold below.
 
  We note first that the following obvious necessary condition holds. 
 
 \begin{lemma}\label{l:ind}
If the matrix $\alpha$ is
irreducible then for any $N>2$ and any $(z_1,v)\in Q_1\cap A_2$ the polynomials $q_j$, $1< j< N$, defined by
\eqref{q_tilde}, are linearly independent.
\end{lemma}

{\it Proof.} Let  us iteratively define the partition $E_i$ of  $\{1,\ldots, M\}$ for\linebreak $i =\{0,\ldots, N\}$ such that $E_0 =\{1,  N\}$ and $j \in E_i$ if $\alpha_{jl} = 0,$ for all $l \in E_{i'}$,
$i' \le i -2$, and there exists $l' \in  E_{i-1}$ such that  $\alpha_{jl'} \ne 0$. Define $E_1$ to be the set of all indices $\{1 < j < N\}$ such that at least one among $\alpha_{j1}$ and $\alpha_{jN}$ is different from zero: since the matrix is irreducible, the set is non-empty. If $(E_0\cup E_1)^c \ne\emptyset$, then there exists $j\in (E_0\cup E_1)^c$ such that $\alpha_{jj'} \ne 0$ for some $j'\in E_1$, otherwise the matrix would be reducible. The procedure continues
until $(\cup_{0\le i\le M}E_i)^c = \emptyset$ and we have that $1\le M \le N -2.$
Consider a linear combination $\sum_{1<j<N}c_jq_j$. By the homogeneity in $(z_2,\ldots,z_{N-1})$ it vanishes identically if and only if
\begin{equation}\label{eq:lin_comb}
    \sum_{1<j<N}c_jz_j\cdot(\alpha_{j1}z_1+\alpha_{jN}v)\equiv 0\,,
    \sum_{1<i,j<N}c_j(\alpha_{ji}-\alpha_{jN}\alpha_{1i})z_j\cdot z_i\equiv 0\,.
 \end{equation}
until $(\cup_{0\le i\le M}E_i)^c = \emptyset$ and we have that $1\le M \le N -2.$
Since $(z_1,v)\in Q_1\cap A_2$  the term in brackets in
the first sum of \eqref{eq:lin_comb} is different from zero for each  $j\in E_1$, so  $c_j=0$ for every $j\in E_1$. Using this
in the second sum of \eqref{eq:lin_comb} we get:
$$
\sum_{n=2}^M\sum_{m=n-1}^M
\sum_{i\in E_m}\sum_{j\in
  E_n}c_j\alpha_{ji}z_j\cdot z_i\equiv0\,.\\ 
$$
This relation holds  if and only if
$
(c_j-c_i)\alpha_{ji}=0
$
 for all $j\in E_n$, $2\le n\le M$, and $i\in E_m$, $n-1\le m\le M$.
We know that  $c_j=0$ if $j\in E_1$.
 Starting from $n=2$ and arguing  by induction in $n$ we find that if $c_i=0$ for all $i\in E_{n-1}$, then $c_j=0$ for all
$j\in E_n$. Indeed,  for any $j\in E_n$ there exists at least one
$i\in E_{n-1}$ such that $\alpha_{ji}\ne0$  by the definition of $E_i$, so relation 
$(c_j-c_i)\alpha_{ji}=0$ implies that $c_j=0$ if $j\in E_n$. That is,  $c_j\equiv0$. 
\qed
\medskip

Further, we have the following necessary condition verified.

For $1 < n < N$ define the intersection $\tilde Q^{(n)} = \cap_{1<j\le n} Q_j$ (and, in particular, $\tilde Q^{(N-1)} =   Q(z_1,v)$). Recall that we continue to suppose that $\al$ is irreducible.

\begin{lemma} $($\cite{Mai}$)$ Reordering  the indices $j \in\{2,\ldots, N-1\} $ one can suppose that the intersection $\tilde Q^{(n+1)}$ is a proper subset both of $\tilde Q^{(n)}$ and $\tilde Q_{n+1}$  for any $1 < n < N $.\end{lemma}

{\em Proof.} We continue to use the partition $E_i$ of  $\{1,\ldots, N\}$ for $i =\{0,\ldots, M\}$ constructed in the proof of Lemma 2.1.

Notice that the equations for $q_j$ , for $1 < j < N,$ can be rewritten as
(we continue to suppose that $z_N=v$)
\be\label{25} q_j = z_j\cdot \left(\sum_{j'\in E_{i-1}} \alpha_{jj'}z_{j'} +\sum_{j'\in E_{i}}
\tilde\alpha_{jj'}z_{j'} +\sum_{j'\in E_{i+1}}
\alpha_{jj'}z_{j'} \right),\; j \in E_i ;
\ee
with $\tilde\alpha_{jj'} = \alpha_{jj'} - \alpha_{jN} \alpha_{1j'}$ (note that $\tilde\alpha_{jj'} =  \alpha_{jj'}$ if $j \in E_i$, for $i > 1$).
Reorder then the indices in  $\{2,  N-1\}$ in such a way that if $j < j'$
and $j \in E_i, j'\in E_{i'}$ , then $i \le i'.$ Fix then $2 \le n < N -1$ and suppose that
$n+1\in E_i.$ The space  $\{z_{n+1} = 0\}$ is contained in $\tilde Q_{n+1}$ but not in $\tilde Q^{(n)},$ since
$ q_2$ is not identically zero ($ 2 \in E_1$), so that $\tilde Q_{n+1}$ is not contained in  $ \tilde Q^{(n)}$. On the other
hand, if $i = 1$ consider the space
$$S_n =\bigcap_{\{2,\ldots,N-1\}\setminus\{n+1\}}
\{z_j = 0\} .
$$
It is contained in  $\tilde Q^{(n)}$, but the restriction of $ q_{n+1}$ to $S_n$  is (note that in this formula and below we sometimes  write $\al_{i,j}$ instead  $\al_{ij} $ by typographical reasons)
$$ q_{n+1}|_{S_n}= z_{n+1} \cdot \left(\alpha_ {n+1,1}z_1 + \alpha_{n+1,N}v + \tilde\alpha_{n+1,n+1}z_{n+1}\right) $$
which is not identically zero, since  $\alpha_ {n+1,1}z_1 + \alpha_{n+1,N}v\ne 0$ ($z_1$ is not proportional to  $v$), so that $S_n \cap \tilde Q^c_{n+1}\ne \emptyset$ and  $\tilde Q^{(n)}$ is not contained in $\tilde Q_{n+1}.$ 

If $i > 1,$ instead, there exists a
chain of indices $j_{i'} \in E_{i'}$ , where $0 < i' < i$ such that   $\alpha_{j_{i'}, j_{i'+1}}\ne 0$,  
$j_i = n + 1 \in E_i,$ and we collect them in the set $I_{n+1} = \{j_{i'}\}_{0<i'<i}$. Indeed, by the very definition of $E_i,$ there exists at least one index 
$l' \in E_{i-1}$ such that $\alpha_{l'{j_i}} \ne 0.$ Fix arbitrarily one of those indices $l'$ and call it $j_{i-1}$. Then, if $i -1 = 1$ we are done, otherwise, since $j_{i-1}\in E_{i-1}$ there exist at least one
$l'' \in  E_{i-2}$ such that $\alpha_{l''{j_{i-1}}}\ne 0: $ fix one of them and continue iteratively until the complete chain is constructed in $i -1 $ steps. Consider the spaces
\be\label{26} V_1 = \{\alpha_{j_1,1}z_1 + \alpha_{j_1,N}v + \tilde\alpha_{j_1j_1}z_{j_1} + \alpha_{j_1j_2}z_{j_2} = 0\} ,\ee 
\be\label{27} V_{i'} = \{\alpha_{j_{i'} j_{i'-1}}z_{i'-1}+ \alpha_{j_{i'} j_{i'+1}}z_{j_{i'+1}} = 0\}\;\; \hbox{\rm for } 1 < i' < i  ,\ee
 and the space
$$S_n =S_n'\bigcap S_n''$$
for
$$
S_n'=\left(\bigcap_{j\in\{2,\ldots,N-1\}\setminus (I_{n+1}\cup \{n+1\}} \{z_j = 0\}\right),\; S_n''=\left(\bigcap_{i'\in I_{n+1}}V_{i'} \right).
$$

  Let $(z_2,\ldots, z_{N-1})$ be a point in  $S_n.$ If
  $j \le n, j \notin I_{n+1},$ then $ q_j = 0,$
since $z_j = 0.$ If $j \notin I_{n+1}$, instead, i.e., $j = j_{i'}$ for some $i' < i,$ note that
$$q_{j_{i'}}|_{S_n'}=z_{j_1}\cdot\left(
\alpha_{j_1,1}z_1 + \alpha_{j_1,N}v + \tilde\alpha_{j_1j_1}z_{j_1}+\alpha_{j_1j_2}z_{j_2}\right),\; \hbox{\rm if } i'=1$$
$$ q_{j_{i'}}|_{S_n'}=z_{j_{i'}}\cdot\left(
\alpha_{j_{i'}j_{i'-1}}z_{j_{i'-1}} + \alpha_{j_{i'}j_{i'+1}}z_{j_{i'+1}}\right),\; \hbox{\rm if } i'\ne 1$$
where the terms in brackets vanish for $\left(z_2, \ldots, z_{N-1}\right) \in V_{i'}.$ Then $S_n$ is
contained in $\tilde Q^{(n)}$, but the restriction of $ q_{n+1}$ to $S_n$  is
$$ q_{n+1}|_{S_n} = \alpha_{j_ij_{i+1}} z_{j_i} z_{j_{i-1}}|_{S''_n}=\pm z_{j_1}\cdot\left(\alpha_{j_1,1}z_1 + \alpha_{j_1,N}v + \tilde\alpha_{j_1j_1}z_{j_1}\right),$$
where $z_{j_1}$ is free. The first equality is a direct consequence of \eqref{25} and  that $z_j = 0$ for $j \notin I_{n+1}\cup \{n + 1\}.$ For the second equality, note out that conditions \eqref{26}-\eqref{27} can be seen as a system of $i-1 $ linear equations in the $i$ unknowns $z_{j_{i'}}$ , for $1 \le i' \le  i.$
Actually, the solutions can be found explicitly considering $z_{j_1}$ as free parameter, by using the fact that $\al_{j_{i'} j_{i'+1}}\in \{\pm 1\}$, so that the conditions $(z_2,\ldots z_{N-1}) \in V_{i'}$ , for $i' > 1,$ become simply $z_{j_{ i'+1}} = -\al_{j_{i'} j_{i'-1}}z_{j_{ i'-1}}/\al_{j_{i'} j_{i'+1}}  = \pm z_{j_{ i'-1}}$. Thus, apart from
the trivial cases $i = 2, 3$, we have that
$$ z_{j_i} = \pm z_{j_{i-2}}= \ldots = \pm\left(\al_{j_1,1}z_1 + \al_{j_1,N}v + \tilde \al_{j_1j_1}  z_{j_1}\right) ,$$
$$z_{j_{i-1}}= \pm z_{j_{i-3}}= \ldots = \pm z_{j_1},
$$
if $i$ is even, and
$$z_{j_{i}}= \pm z_{j_{i-2}}= \ldots = \pm z_{j_1},
$$
$$ z_{j_{i-1}} = \pm z_{j_{i-3}}= \ldots = \pm\left(\al_{j_1,1}z_1 + \al_{j_1,N}v + \tilde \al_{j_1j_1}  z_{j_1}\right) ,$$
if $i$ is odd. Thus $S_n \bigcap \tilde Q^c_{n+1}
 \ne \emptyset$ and $\tilde Q^{(n)}$ is not contained in $\tilde Q_{n+1} ,$ which completes the proof.
 Note also that since $S_n$ belongs to the intersection
of all quadrics different from $\tilde Q_{n+1},$   we get that for any $1 < j < N,\; \bigcap_{k\ne j} \tilde Q_{k}$ is not contained in $\tilde Q_{j+1}.$ \qed\medskip

Unfortunately, that necessary condition becomes sufficient only if all intersections $\tilde Q^{(n)}$  are irreducible which is not clear.\smallskip
 
Finally, let us note that the conjecture is true in small dimensions:
\begin{proposition}\label{01-34}
  Let  $N\in\{3,4\}$. Then  MC holds true.
\end{proposition}

{\em Proof}.
{$\bf   N=3. $}  
If $N=3$ then $N-2=1$ and we have only one non-trivial equation for a fixed  $(z_1,v)\in Q_1\cap A_2$ and the conclusion is immediate. 

{ $\bf   N=4.$} We have to show that the codimension of  intersection of the two quadrics is two 
 (and not one). We begin with the next evident  lemma.

  {\begin{lemma}\label{l}
 Let $\mathcal{Q}_1=\{ q_1=0\},\mathcal{Q}_2=\{ q_2=0\}$ be two  linearly independent quadrics over $K$.   Then the codimension of $\mathcal{Q}_1\cap \mathcal{Q}_2$ is one if and only if  $ q_1$ and $ q_2$   have a mutual affine linear factor $l(x)$.\qed \end{lemma}

Now to prove Proposition \ref{01-34} we have to consider  two cases. If one of the polynomials $ q_2$ and $q_3$ is linear,  the codimension is two since  they are linearly independent. 

In the second case both $q_2$ and $q_3$ are quadratic. Then  the codimension is still two since the polynomials $q_{j} (z_2,\ldots,z_{N-1};z_1,v), \, j=2,3,$     are irreducible. Indeed, they can be written 
$q_{j} =z_j\cdot l_j(z_2,\ldots, z_{N-1};z_1,v)$ , 
 where $l_j$ are surjective affine functions  $l_j:K^{d(N-2)}\longrightarrow K^d$. But such scalar product cannot vanish for $d\ge 2>1$ on a hyperplane $H \subset K^{d(N-2)} $ which  would be the case for a reducible quadric. Indeed we have either:

 a) the coefficient $a$ of $z_j$ in $l_j$ is non-zero, or

b) it is zero but then the coefficient $b$ of some other  $z_i$  is non-zero. 

 In case a) take the 2-dimensional plane $P(x_1,x_2)$ in
  the whole space  generated by some two  vectors from the
  $z_j$-space. In particular, we can choose the first basic vector
  parallel to $\alpha_{1j}z_1+\alpha_{Nj}v\neq 0$ (otherwise case a)  is impossible), and the restriction of $q _{j}=0$ on $P$ is then
  $a(x_1^2+x_2^2)+c_1x_1=0$, with $c_1\neq 0$, which is isomorphic  to $x_1^2+x_2^2=C\neq 0$, and 
this   quadric in $P(x_1,x_2)$   cannot contain $P(x_1,x_2)\cap H$ (a line or the whole  $P(x_1,x_2)$).

Similarly, in case b) we take the 4-dimensional vector subspace $P'$  generated by the two first basic vectors in the $z_j$ space and  the  two first basic 
vectors in the $z_i$ space. The restriction of $q _{j}=0$ on $P'$ is then $b(x_1y_1+x_2y_2)+c_1x_1+c_2x_2=0$,  isomorphic to $x_1y_1+x_2y_2=C$
 which   can not contain $P'(x_1,x_2,y_1,y_2)\cap H$.  
This finishes the proof for $N=4$.

The above proof implies also the following result

\begin{lemma}\label{l:irr}
The polynomials $q_{j},\;1< j< N$, defined by
\eqref{q_tilde}, are irreducible for an irreducible  matrix $\alpha$.
\end{lemma}

{\em Remarks.} \smallskip

 I. Arguments, similar to the case $N=4$ permit to deal also with the cases $N=5,6$  confirming the conjecture for them, which needs considering numerous various cases for the matrix $\al.$ However, the combinatorics becomes complicated for larger $N$, and for $N\ge 7$ this method seems to be  unfeasible.

 II. For the purposes of the Zakharov--L'vov theory, as a simple analysis   of the arguments below shows, it is sufficient to prove MC for $\K=\C.$
\newpage
\section{Arithmetic 
consequences} \lbl{sec:arithm} In this section we deduce from \eqref{mc} some arithmetic consequences concerning the asymptotic behavior of certain interesting  sums.
 \subsection{Main sums and their estimates}\lbl{sec:numbertheory}
Let us now define those sums and formulate some   estimates for them.

We fix two integers $N\ge 2, d\ge 3$ and consider a function $\Phi:
\mathbb R^{Nd}\to \mathbb \R$ which is sufficiently smooth and sufficiently fast decaying at infinity (see below for detailed assumptions).
Our goal is to study asymptotic as $L\to\infty$ behaviour of the sum 
 \be\lbl{S_L}
 S_{L,N}(\Phi) := L^{N(1-d)} \sum_{z: \  \omega_j (z) =0 \, \forall j} \Phi(z), \ee
 where $\omega_j$ are defined by \eqref{om} over $\K=\C.$
For a function $f\in C^k(\mathbb R^m)$ and
$n_1,n_2\in \mathbb N$, $n_1\le k$ we define
$$\|f\|_{n_1,n_2}= \sup_{z\in \mathbb R^m} \max_{|\alpha|\le n_1} |\p^\alpha f(z)| \lan z\ran^{n_2}\,,\;\;\lan z\ran:=\max\{1,|x|\}.
$$

The first crucial result concerns the case $N=2$. Note that for $N=2$ the set $\{z:  \omega_j (z) =0 \, \forall j\}$ over which we take the summation in \eqref{S_L} takes the from $\{z\in \mathbb
	Z^{2d}_L: \; z_1\cdot z_2=0\}$ since $\om_1(z)=-\om_2(z)=\al_{12}z_1\cdot z_2$ where $\al_{12}\ne 0$.

Now, the asymptotic for the sum $S_{L,2}(\Phi)$ immediately follows from Theorem~1.3 in \cite{nt} with $\eps=1/2$.
\begin{theorem}\lbl{t:numbertheory}
Let    $N_1(d):=4d(4d^2+2d-1)$, $N_2(d):=N_1+6d+4$. If
$\|\Phi\|_{N_1,N_2}< \infty$, there exist constants
$C_d\in (1,4/3)$, $K_d>0$ such that
  \be\label{H-B}
  \left|S_{L,2}(\Phi) -  C_d \int_{\Sigma_0}
{\Phi(z)} \,   
d\mu^{\Sigma_0}  
\right|\le K_d {\|\Phi\|_{N_1,N_2}}{L^{5/2-d}} \,,
  \ee
where  $\Sigma_0$ is the quadric  $ \{z\in\R^{2d}: z_1 \cdot
z_2=0\}$ equipped with the measure  
$d\mu^{\Sigma_0} =  \left(dz_1 dz_2/\sqrt{|z_1|^2+|z_2|^2}\right) \!\mid_{\Sigma_0}. $  
\end{theorem}

The integral in \eqref{H-B} converges if $\Phi(z)$ decays at infinity:
\be\label{int_est}
\Big|\int_{\Sigma_0} {\Phi(z)} \, \mu^{\Sigma_0} (dz_1dz_2) \Big| \le C_r \| \Phi\|_{0,r} \quad\text{if} \;\; r>2d-1,
\ee
see Proposition 3.5 in \cite{DK1}. 

Assuming \eqref{mc} we can     deduce from Theorem \ref{t:numbertheory}   another result, whose proof is given in  to  Section~\ref{sec:intersection} below :
\begin{theorem}\lbl{t:countingterms} Let \eqref{mc} hold for a given $N\ge 3$ and ${\mathbb K}={\bar{\mathbb F}_p},$ the algebraic closure of a finite field.
  Then there exists a constant $C_{d,N}$ such that
  $$
  \big| S_{L,N}(\Phi) \big| \le C_{d,N}\|\Phi\|_{0,\bar N}\,,
  $$
  for $\bar N:=\lfloor N/2 \rfloor N_2(d)+(N-2)(d-1)+ 1$ and $N_2$
  defined in Theorem~\ref{t:numbertheory}. 
\end{theorem}
Note that this is, in particular, the case for $N=3$ and $4$. Below in this section we  suppose that  \eqref{mc} holds for our  $N\ge 3$, ${\mathbb K}={\bar{\mathbb F}_p}$ and any prime $p$.
 
 The  proof of Theorem \ref{t:countingterms} is based on the following counting result
 \begin{lemma}\label{l:intersection} Let $R>0,$ and let $\;Q_{1,L}:=Q_1\cap \Z^{2d}_L,\;\tilde Q_{L}:=\tilde Q\cap \Z^{(N-2)d}_L.\;$ \linebreak
Assume that the matrix $\al$ is irreducible. Then for any 
$(z_1,v)\in Q_{1,L}\cap A_2$, the cardinality $s(R,\tilde Q,L)=|S(R,\tilde Q,L)|$ of 
$$S(R,\tilde Q,L) = \tilde Q_L(z_1,v)  \cap B_R^{(N-2)d}$$  verifies
$$
s(R,\tilde Q,L)\le (2RL)^{(N-2)(d-1)}\,,
$$
$B_R^{(N-2)d}$ being the cube $\big\{|z|_\infty< R\big\} \subset{\mathbb \R}^{(N-2)d}.$
\end{lemma}
 Thus, in  Lemma~\ref{l:intersection} we want to estimate the number of integer points on a quadrics' intersection inside a large box. The idea is to embed the (integer points of the) box in an affine space over a large finite field and then apply  finite fields' version of the  Bezout theorem to estimate the cardinality.   
 \subsection{Finite fields' Bezout theorem}
The   {\em Bezout theorem} states in its  most elementary version that for an affine algebraic set (AAS below) $X$ one has
$$\deg X\le \Pi_{i=1}^s \deg F_i,$$
where $X=\bigcap_{i=1}^s H_i $ for hypersurfaces $H_i=\{F_i=0\}$ defined by forms $F_i.$ We use it's finite fields'  version  \cite[Corollary 2.2]{LR}.

\begin{theorem}\label{bz} Let ${\mathbb K}=\bar{\mathbb F}_p$ and let $X=\bigcap_{i=1}^s H_i $ be an AAS   with $\dim X=r$ and $d_i=\deg F_i, i=1,\ldots, s$. Then
$$ |X\cap{\mathbb F}_p |\le p^r \Pi_{i=1}^s d_i.$$
\end{theorem}

\subsection{Proof of Lemma~\ref{l:intersection}}
Let $q_1,\ldots, q_s$, $s\ge 1$, be  polynomials of degree at most two in $m\ge s$ variables,  $q_i\in \Z[X_1,\ldots, X_m]$,
with $q_i(0)=0,$\linebreak $ i=1,...,s$. 
Consider the geometric quadrics $Q_i=\{x\in \R^m: q_i(x)=0\}$ and their intersection $Q=\cap_{i=1}^s Q_i$. The latter is not empty since $\{x=0\}\in Q$. 
 
Let
$B_M^{m}\subset \mathbb R^{m}$ 
be the open cube $\left|x\right|_\infty<M$, with $M\ge 1$.  Consider the set
$$
S_m(M,Q) =Q\cap \Z^{m}\cap B_{M}^m .
$$
 
Let $p$ be a prime and  $q^{(p)}_{i}(X)\in {\mathbb F}_p[X_1,\ldots, X_m]$ denote the polynomials $q_{i}(X)\mod \,p$ over the finite field ${\mathbb F}_p$. 
Consider the sets (recall that $K=\bar{\mathbb F}_p $ is the algebraic closure of ${\mathbb F}_p$)
$$Q^{(p)}_{i}=\{x\in  K^m: q^{(p)}_{i}(x)=0\}$$
and their intersection $Q^{(p)}=\cap_{i=1}^s Q^{(p)}_{i}$. 
We will be interested mainly in  the cardinality of $Q^{(p)}({\mathbb F}_p):=Q^{(p)}\cap{\mathbb F}_p^m$
as a tool to estimate $|S_m(M,Q)| $.

\begin{proposition}\label{01} Let $M\ge 1$ be large and  suppose that for a  prime\linebreak $p\in [2M+1, 2M +o(M)] $,  one has $$ \quad \quad\quad\quad \quad \dim \,Q^{(p)}=m-s\hskip 3 cm $$ (i.e.,  the quadrics $Q^{(p)}_{i}$ intersect properly).
 Then
  $$|S_m(M,Q)|  \le   \big(2^{m}+o(1)\big)M^{m-s} \, . $$
  \end{proposition}
 {\em Proof.}  Let
 $ \Pi: S_m(M,Q)\longrightarrow {\mathbb F}_p^m $
be defined by
$$\Pi(x_1, \ldots, x_m)=(x_1 \!\!\mod p, \ldots, x_m\!\!\mod p).$$
Then $ \Pi$ is injective  and its image is contained in $Q^{(p)}\cap{\mathbb F}_p^m\subset {\mathbb F}_p^m$. Indeed, the last assertion is clear and 
the injectivity is established as  follows: if  
\begin{equation*}\hskip1 cm (x'_1\!\mod p, \ldots, x'_m \!\mod p) =(x_1\!\mod p, \ldots, x_m\!\mod p)\hskip1 cm \end{equation*}
 but $x'\ne x$, then for some $  i\in\{1,...,m\}$ we have  $x'_i\!\mod p =x_i\!\mod p $, but $\,x_i'\ne x_i$. 
 Consequently, $|x_i-x'_i|\ge p>2M$ which contradicts the condition $x_i,x'_i\in B_{M}^m$.
Applying then Theorem \ref{bz} to $X=Q^{(p)}$ over ${\mathbb K}=\bar{\mathbb F}_p$ we get the conclusion since
$$|S_m(M,Q)|\le |Q^{(p)}({\mathbb F}_p)|\le
    2^sp^{m-s}=\big(2^{m}+o(1)\big)M^{m-s}\,.
$$
\qed

\noindent 
Such primes $p$   exist for any large $M$ by the Prime number theorem.\smallskip

Now we pass to the proof of Lemma~\ref{l:intersection} and set
$$
s(R,{\tilde { Q_L}},L)=\big|  \tilde Q_L(z_1,v)  \cap B_R^{(N-2)d}\big|\, .
$$
Consider the set
\begin{equation*}
\begin{split}
\quad S'\big(R,{\tilde {\mathcal Q}},L\big) &=
{\tilde {\mathcal Q}}\cap \Z^{(N-2)d}\cap B^{(N-2)d}_{RL}\,,\qquad {\tilde {\mathcal Q}}=\tilde Q_L(Lz_1,Lv), 
\end{split}
\end{equation*}
and let $s'(R,{\tilde {\mathcal Q}},L)$ be its cardinality. Then 
 $s' (R,{\tilde {\mathcal Q}},L) = s(R,{\tilde { Q}},L)\,,$
since  the map 
$(z_2,\ldots,z_{N-1}) \mapsto (L z_2,\ldots,L z_{N-1})$ is  a
bijection between  $S(R,{\tilde {\mathcal Q}},L)$ and $ S'(R,{\tilde {\mathcal Q}},L)$.

Thus, we can estimate $s(R,{\tilde {\mathcal Q}},L)$ through Proposition~\ref{01} with  $M=RL$ and  $m=(N-2)d, s=N-2$ to deduce a proof of Lemma~\ref{l:intersection}   since our main conjecture \eqref{mc} implies that $\dim \,Q^{(p)}=m-s$. 
 \subsection{Proof of Theorem~\ref{t:countingterms}}\label{sec:intersection}
Let us  define $Q_j:=\{z\in \R^{Nd}:\, \omega_j(z)=0\}$ and consider the  intersection  $Q=\cap_{j=1}^N Q_j$. 
  Note that $Q=\cap_{j=1}^{N-1} Q_j$ since 
  $\om_1+\ldots+\om_N=0$.
 Denote by $B_R^{Nd}$ the open cube $|z|_\infty<R$ in $\R^{Nd}$.

\begin{proposition}\label{p:finite_support}
If for $w:\R^{Nd}\longrightarrow \R$ one has 
$\supp(w)\subset B_R^{Nd},$ and  $ \;\|w\|_\infty< \infty,$  then
$$ \left|\sum_{z\in Q\cap \cZ} w(z)\right| \le 
C(N,d) R^{\lfloor N/2\rfloor N_2(d)+(N-2)(d-1)}L^{N(d-1)}\|w\|_\infty\,
$$ 
for any $R\ge1$ where $N_2$ is defined in Theorem~\ref{t:numbertheory}.
\end{proposition}

{\it Proof.}
Recall that $ Q_1=\{(z_1,v)\in\mathbb R^{2d}:z_1\cdot v=0\},\;Q_{1,L}=Q_1\cap \Z^{2d}_L$, $\tilde Q_{L}=\tilde Q\cap \Z^{(N-2)d}_L$, and
 $A_2=\{(z_1,v)\in\mathbb R^{2d}: z_1\neq 0, v\neq 0\}, $  $\cZ$ being the set with no zero coordinate
Then
\begin{equation}\label{eq:sum_intersec}
  \left|\sum_{z\in \cZ\cap Q } w(z)\right|\le C(N,d)\|w\|_\infty \sum_{(z_1,v)} 1  \times\sup_{(z_1,v)\in Q_{1,L}\cap A_2}\;
\sum_{(z_2,\ldots,z_{N-1}) }1\,.
  \end{equation}
where in the first sum  $(z_1,v)$ runs over  $Q_{1,L}
  \cap B_{ (N-1)R}^{2d} $ while in the second sum
  $(z_2,\ldots,z_{N-1})$ runs over $\tilde Q_L(z_1,v) \cap B_R^{(N-2)d}. $\smallskip

 To estimate the first  sum, we take any smooth function $w_0(x) \ge0$, equal one for $x\le1$ and vanishing for
$x\ge2$. Then 
$$
\sum_{(z_1,v)\in Q_{1,L}\cap B_{ (N-1)R}^{2d}} 
1 \le \sum_{(z_1,v)\in Q_{1,L}} w_R(z_1,v),
$$
where $w_R(z_1, v) :=  w_0\Big( |(z_1,v)| / \big(  (N-1) R\sqrt{2d} \Big)$.  Since for $R\ge1$ we have 
$\| w_R\|_{N_1, N_2} \le CR^{N_2}, $
and by Theorem~\ref{t:numbertheory} and \eqref{int_est} get 
 $$ \sum_{(z_1,v)\in Q_{1,L}\cap B_{ (N-1)R}^{2d}} 1 \le 
C  L^{2(d-1)} \big[ R^{2d} + R^{N_2} L^{-d+5/2}\big] \le 
C'  L^{2(d-1)}  R^{N_2},
 $$  
where $C, C'$ depend on $d, N, N_1$ and $ N_2$. 

For the second sum(s) of \eqref{eq:sum_intersec} we  use 
 Lemma~\ref{l:intersection}.\smallskip

This completes the proof of Proposition~\ref{p:finite_support} in the case of irreducible matrix $\al$: indeed, applying MC we get 
\be\lbl{sum_w_est}
\left|\sum_{z\in \cZ\cap Q } w(z)\right|
\le C(N,d)\|w\|_\infty R^{N_2 + (N-2)(d-1)} L^{N(d-1)}.
\ee
If the matrix $\al$ is not irreducible, it can be reduced through
permutations to a block diagonal matrix with $m$
blocks which are irreducible square matrices of sizes $N_i$ satisfying $\sum_i N_i=N$. Since $N_i\ge 2$ (otherwise there would be a zero row or column), $ m \leq \lfloor N/2 \rfloor$.
Applying  \eqref{sum_w_est} to each block we get assertion of the proposition.
\qed\smallskip

Let $\varphi_0(t) = \chi_{(-\infty,1]}(t)$ and  $\varphi_k(t) = \chi_{(2^{k-1},2^k]}(t)$ for $k\ge 1$. Then 

$1=\sum_k\varphi_k(t)$ and $\Phi=\sum_{k=0}^\infty w_k(z)\,,\quad w_k(z)=\varphi_k (|z|_\infty|)\Phi(z)\,.$
Since $\mathrm{supp}\,w_k\subset B_k=\{|z|_\infty\le 2^k\}$ and $\|w_k\|_{\infty}\le
C2^{-k\bar N}\|\Phi\|_{0,\bar N}$, one gets by Proposition~\ref{p:finite_support} that
$$
|S_{L,N}(\Phi)|\le C(N,d)\|\Phi\|_{0,\bar N} \sum_{k=0}^\infty
2^{k(\lfloor N/2 \rfloor N_2+(N-2)(d-1)-\bar N)}\,, 
$$
which converges if $\bar N>\lfloor N/2 \rfloor N_2+(N-2)(d-1)$. This
completes the proof of Theorem~\ref{t:countingterms}.
\qed

\section{Zakharov-L'vov model for wave turbulence.}
 Recall then some basic facts on the Zakharov-L'vov stochastic model for wave turbulence (WT) \cite{ZLF,DK1,DK2}. For simplicity of exposition we suppose that $d\ge 3$ until the end of our paper.

\subsection{Classical setting.}
Let ${\mathbb{T}}^d_L={  \mathbb{R}   }^d/(L{  \mathbb{Z}   }^d)$ be the $d$-dimensional torus, 
of  period $L\geq 2$. One sets   $\|u\|^2 =L^{-d}\int_{{{\mathbb T}}^d_L} |u(x)|^2\,dx\,,$ and  writes the Fourier series  
\begin{equation}\lbl{Fourier-def}
u(x)= L^{-d/2} \,{\sum}_{s\in{{\mathbb Z}}^d_L} v_s e^{2\pi  i s\cdot x}, 
\qquad\qu{{\mathbb Z}}^d_L = L^{-1} {{\mathbb Z}}^d\,.
\end{equation}
Therefore,
$$
v=\cF(u), \quad v_s = L^{-d/2}\,\int_{\T^d_L} u(x) e^{-2\pi i s\cdot x}\,dx \quad \text{for}\; s\in \Z_L^d, 
$$
so the  Parseval identity reads
$\|u\|^2 = L^{-d}\,{\sum}_{s\in\Z^d_L} |v_s|^2.$

Consider  the cubic NLS equation with modified non-linearity 
\be\label{NLS}
\frac{{{\partial}}}{{{\partial}} t}  u +i\Delta u- i\la \,\big( |u|^2 -2\|u\|^2\big)u=0 , 
\qquad x\in {\mathbb{T}}^d_L,
\ee
where $u=u(t,x)$,  $ \Delta=(2\pi)^{-2}\sum_{j=1}^d ({{\partial}}^2 / {{\partial}} x_j^2)$ and 
$\la\in(0,1]$ is a small parameter. One  studies  solutions $u(t,x)$ with $\|u(t,\cdot)\|\sim 1$ as $L\to\infty$. The modification of the non-linearity by the term $2i\la\|u\|^2u$  keeps the main features of the standard cubic NLS equation, reducing  some non-crucial technicalities; see  the introduction to \cite{DK1}. 

The objective of WT is to study solutions of \eqref{NLS}  under the limit $L\to\infty$ and $\la\to 0$ on long time intervals;  some references  containing    different (but consistent) approaches may be found in \cite{ZLF, Naz11}.
\subsection{Zakharov-L'vov setting.} When  studying  eq. \eqref{NLS}, crucial  mechanism, used by WT physicists,  is "pumping  energy to low modes and dissipating it in high modes".  
To make this rigorous, Zakharov and L'vov    proposed  to
 consider the NLS equation \eqref{NLS} dumped by a (hyper) viscosity and driven by a random force:
\begin{equation}\label{ku3s}
\frac{\partial u}{\partial t} +i  \Delta u  - i\la\,\big( |u|^2 -2\|u\|^2\big)u =  -\nu\frak A(u)    + \sqrt\nu\frac{{{\partial}}}{{{\partial}} t} \eta^\omega(t,x).
\end{equation}
Here $\nu\in(0,1/2]$ is another  small parameter, which should
agree properly with $\lambda$ and $L$ for certain linear {\em dissipation  operator} $\frak A,$ 
(we do not give here it's exact form) and the random noise  $\eta^\omega$  is given by a  Fourier series
$$
\eta^\omega(t,x)=
L^{-d/2} \, {\sum}_{s\in\Z^d_L} b(s) \beta^\omega_s(t) e^{2\pi  i s\cdot x},
$$
where
$\{\beta_s(t), s\in{{\mathbb Z}}^d_L\}$ are standard independent complex Wiener processes and	 $b(x)$  is a  Schwartz  function on ${{\mathbb R}}^d\supset \Z^d_L$. 

Solutions $u(\tau)$ of \eqref{ku3s} are random processes in
the space $H=   L_2(\T^d_L, \C), $
equipped with the norm $ \|\cdot\| $.
If the coefficients of $\frak A$ are sufficiently smooth, equation \eqref{ku3s} is known to be well posed and  $\EE\|u(\tau)\|^2$ is bounded uniformly in $\tau$ and $L,\nu,\la$, once $\EE\|u(0)\|^2$ is bounded uniformly in these parameters.

It is convenient to pass    to the slow time $\tau=\nu t$ and write  eq.~\eqref{ku3s} as 
\begin{equation}\label{ku3}
\begin{split}
\dot u +i \nu^{-1} \Delta u  - i\rho\,\big( |u|^2 -2\|u\|^2\big)u &= -\frak A(u)    + \dot \eta^\omega(\tau,x),\\
\eta^\omega(\tau,x)&=
L^{-d/2} \, {\sum}_{s\in\Z^d_L} b(s) \beta^\omega_s(\tau) e^{2\pi  i s\cdot x}\,.
\end{split}
\end{equation}
Here $\rho=\la\nu^{-1}$, the upper-dot stands for $d/d\tau$  and $\{\beta_s(\tau), \, s\in\Z^d_L\}$ is another set of standard independent complex Wiener processes;  $\rho$, $\nu$ and $L$ are used 	 as  parameters of the equation.

In the context of equation \eqref{ku3}, 
the objective of WT is to study its solutions $u(\tau)$ when
$L\to\infty \;\;\; \text{and}\;\; \; \nu\to 0.$ One begins (cf.  \cite{DK1,DK2}) by a  study of formal decompositions in $\rho$ of solutions to   \eqref{ku3} and of their  energy spectra $N_s$ (see the  definition below in Section 5.1) under that limit. \smallskip

Then  one should choose a relation between the  parameters $\nu$ and $L,$ since the theory postulates no such relation.  Using the  choice of \cite{seclim}, namely, 
    {\em firstly} $\; \nu\to 0, $ {\em and then } $ L\to\infty$,
  the main result   is that the behaviour of  principal terms of the decomposition in $\rho$  for the energy spectrum $N_s$  (the  scaling  being $\rho\sim L $) is governed by a {\it modified wave kinetic equation} (mWKE, see below Section  5.3), which is  similar to (but different from) the WKE arising in physics.
  
\section{ Applications to the large period limit.}
In this section we suppose that   \eqref{mc} holds and give some consequence for WT. Thus, the results of Section 3 are valid. Eventual proofs of the results  below can be   rather  close  to those in \cite{seclim} where the cases of $M=2$ and $3$ are treated using $N=3$ and $4$ in Proposition 2.3.  
\subsection{The limit of discrete turbulence}
One first considers the limit 
\be\lbl{lim_dt} 
\nu\to 0 \qmb{while $L$ and $\rho$ stay fixed.}
\ee
It is known as
\textit{the limit of discrete turbulence} (see \cite[Section 10]{Naz11}). Taking the Fourier transform of  ~\eqref{ku3},  and passing to the so-called {\it interaction representation}, 
$$\ai_s(\tau)=v_s(\tau) e^{-i\nu^{-1}\tau |s|^2}, \qu s\in\Z^d_L,$$
we get the equation
\begin{equation}\label{ku44}
 \dot \ai_s + \gamma_s \ai_s= i\rho Y_s(\ai,\nu^{-1}\tau) +b(s) \dot\beta_s,  \qquad s\in{{\mathbb Z}}^d_L,
 \end{equation}
where
$$Y_s(\ai, t )=L^{-d} \Big({\sum}_{1,2,3} \dep \ai_{1}\ai_{2} \bar \ai_{3}
e^{  i t \oms} - |\ai_s|^2\ai_s\Big), $$
$$
\oms:= |s_1|^2+|s_2|^2-|s_3|^2 - |s|^2=2(s_1-s)\cdot (s_2-s),
$$
$\{\beta_s\}$ being yet   another set of standard independent complex Wiener processes ($\ai_j$ standing for $\ai_{s_j}$). Note that the energy spectra  of $v_s(\tau)$ and $\ai_s(\tau)$ coincide: 
$$ N_s(\tau)=\EE|v_s(\tau)|^2=\EE|\ai_s(\tau)|^2. $$
The limiting dynamics in   \eqref{ku44} under  \eqref{lim_dt} 
is governed by  the \textit{effective equation} of discrete turbulence. The latter  has the  form  \eqref{ku44} with the modified nonlinearity  $Y^{res}$, in which  the sum is taken only over resonant vectors $s_1,s_2,s_3$:
\begin{equation}\label{ku5}
\begin{split}
\dot \ai_s+\gamma_s \ai_s &= i\rho Y_s^{res}(\ai)
+b(s) \dot\beta_s
\,,\quad s\in{{\mathbb Z}}^d_L  \,,\\
Y_s^{res}(\ai)&=L^{-d}\Big( {\sum}_{1,2,3} \dess \delta(\omega^{12}_{3s}) 
\ai_{1} \ai_{2} \bar \ai_{3} - |\ai_s|^2\ai_s\Big)\,.
\end{split}
\end{equation}
Here   $\delta(\omega^{12}_{3s})=1$ if $\omega^{12}_{3s} = 0$ and $\delta(\omega^{12}_{3s})=0$ otherwise.
\subsection{Quasisolutions}
 In  WT it is traditional to analyse the quasisolution (defined below) instead of the solution itself, postulating that the former well approximates the latter; see Introduction in \cite{DK1} for a discussion.\smallskip

{\it Quasisolutions and their energy spectra.} Assume for simplicity that initially the system was  at rest, 
\be\label{in_cond}
\ai_s(0)=0 \qquad \forall s\in\Z^d_L.
\ee
and decompose  formally the corresponding solution of \eqref{ku5}  in $\rho$,
\be\label{decompor}
\ai(\tau)=\ai^{(0)}(\tau)+\rho \ai^{(1)}(\tau)+ \rho^2 \ai^{(2)}(\tau) + \dots \,,\quad\qquad 
\ee
$\ai^{(k)}(\tau)=\ai^{(k)}(\tau;L),\,\ai^{(k)}(0)=0.$
Then $\ai^{(0)}(\tau)$
satisfies the linear equation \linebreak  $
\dot \ai_s^{(0)} + \gamma_s \ai_s^{(0)}
= b(s) \dot\beta_s  \,,
 $
for $s\in{{\mathbb Z}}^d_L$  and thus (by Duhamel)
$$
\ai^{(0)}_s(\tau) = b(s) \int_{0}^\tau e^{-\gamma_s(\tau-l)}d\beta_s(l),
$$
while $\ai^{(1)}$ satisfies $\dot \ai^{(1)}_s + \gamma_s \ai^{(1)}_s 
=  i Y_s^{res}(\ai^{(0)})\,. $
Therefore
 $$
\ai^{(1)}_s(\tau) =\frac{i}{L^d}  \int_{0}^\tau  e^{\ga_s(l-\tau)}
\left(\sum_{1,2,3}\dess \delta(\omega^{12}_{3s}) (\ai^{(0)}_1
\ai^{(0)}_2{\bar \ai}^{(0)}_3) -|\ai_s^{(0)}|^2\ai_s^{(0)} \right)\!(l)dl\,
$$
is a Wiener chaos of third order (see \cite{Jan}). 
Similarly for $ n\ge1$  
$$ 
	\ai^{(n)}_s (\tau) 
	= \frac{i}{L^d} \sum  \int_{0}^\tau  e^{\ga_s(l-\tau)}  
	    \left( \sum_{1,2,3}  \dess \delta(\omega^{12}_{3s})
	\ai_{1}^{(n_1)} \ai_{2}^{(n_2)} {\bar \ai_{3}}^{(n_3)}  - \ai_s^{(n_1)} \ai_s^{(n_2)} \bar\ai_s^{(n_3)} \right)\,(l)dl, 
 $$
 (the outer summation is taken over  ${n_1+n_2+n_3=n-1}$) is a Wiener chaos of order $2n+1$. 	

Next one considers the $M$-order truncation of the series \eqref{decompor},
$$\mathcal{A}_s(\tau; L)= 
\mathcal{A}_s(\tau) = \ai^{(0)}_s(\tau) +\rho \ai^{(1)}_s(\tau) +\rho^2 \ai^{(2)}_s(\tau)+\ldots+\rho^M \ai^{(M)}_s(\tau)\,,
 $$
which we call the ($M$-order) {\it quasisolution}  of     \eqref{ku5}, \eqref{in_cond}. 

The main goal here   is to study the behavior of the energy spectrum 
  $\mathfrak{n}_{s,L}(\tau) = \EE |\cA_s(\tau;L)|^2 $ 
of  $\mathcal{A}(\tau)$,  as $L\to\infty$, and to show that it  has a non-trivial asymptotic behavior for  $\rho\sim L$.

  Therefore, one  assumes that  $\rho=\eps L, $ 
where $0<\eps\leq 1/2$ is a small but fixed constant. 
Then the energy spectrum $\no_{s,L}$ for $s\in\Z_L^d$  expands as 
$$\no_{s,L}(\tau)=\no_{s,L}^{(0)}(\tau) + \eps \, \no_{s,L}^{(1)}(\tau)+\eps^2 \no_{s,L}^{(2)}(\tau) +\eps^3 \no_{s,L}^{(3)}(\tau) +\ldots +\eps^{M} \no_{s,L}^{(M)}(\tau), 
$$ 
where  
 $$\no_{s,L}^{(k)}(\tau) =   L^k \sum_{\substack{k_1+k_2=k\leq M }}   \ai^{(k_1)}_s(\tau) \bar \ai^{(k_2)}_s(\tau)\, $$
 for any $k\ge 0$.
In particular,   $
\no_{s,L}^{(0)} (\tau)= \EE |\ai_s^{(0)}(\tau)|^2
= {b(s)^2}\big( 1-e^{-2\ga_s\tau}\big)/{\ga_s} ,$
and a simple computation shows that $\no_{s,L}^{(1)}(\tau)\equiv 0$.
For higher order terms one proves that 
 $ \no_{s,L}^{(2)}\sim 1\qnd |\no_{s,L}^{(k)}|\,\lesssim 1 \quad\mbox{as }L\to\infty \mbox{ uniformly in }\tau\geq 0$
for any  $3\le k\le M$. Thus, the parameter $\eps$ measures the 
properly scaled amplitude of the solutions, and it should be small for the methodology of WT to apply. The  term $\eps^2\no_{s,L}^{(2)}$ is the crucial non-trivial component of the energy spectrum $\no_{s,L}$ while the terms $\eps^k\no_{s,L}^{(k)},$ are perturbative, which  agrees with the prediction of physical works concerning WT.
\subsection{(Modified) wave kinetic equation (mWKE)} Therefore,  to study the limiting as $L\to\infty$ behavior of the energy spectrum $\no_{s,L}(\tau)$ up to an error of size  $\eps^{2M}$ it remains to describe the behavior of its principal component\linebreak 
$\no^{(0)}_{s,L}(\tau)+\eps^2\no_{s,L}^{(2)}(\tau)$. 
In fact, the latter is governed by a mWKE; to state it one considers the {\it resonant quadric}
$$ \Sigma_s=\{(s_1, s_2)\in\R^{2d}: \, (s_1-s)\cdot(s_2-s)=0\},
$$
 equipped by the measure $d\mu: =    \big(ds_1 ds_2/\sqrt{|s_1-s|^2+|s_2-s|^2} \big)\mid_{\Sigma_s}.$ 

 Then our damped/driven non-autonomous mWKE reads as
\be\label{kin eq}
\dot \zz_s (\tau) = -2\ga_s\zz_s +\eps^2 K_s(\tau)(\zz_s) +2 b(s)^2, \qquad  \;\; \zz_s(0)=0,
\ee
where $\tau\ge0$ and $s\in\R^d$. The non-autonomous cubic {\it wave kinetic } operator $K(\tau)$, acts as follows  for any $\tau\geq 0:$ 

  $ K_s (\tau)y_s={  4}C_d\int_{\Sigma_s} d\mu\Bigl(\cZ^4 y_1y_2y_3 +\cZ^3y_1y_2 y_4- \cZ^2 y_1 y_3y_4-  \cZ^1 y_2y_3y_4\Bigr).$ 

Here $y_j:=y_{s_j}$ with  $s_4:=s$ and  $s_3:=s_1+s_2-s$, $C_d$ being a constant, and   
$\cZ^j=\cZ^j(\tau_0)=\cZ^j(\tauz;\vs)$   are given for a fixed  $\tauz> 0$ and  $\,m,j=1,2,3,4$   by    
 
  $ \cZ^j(\tauz;\vs) := \int_0^{\tauz}dl\,\e^{\ga_{s_j}(l-\tauz)} \prod_{m\neq j} \big(\sinh(\ga_{s_m}  l)/{\sinh(\ga_{s_m} \tauz)}\big), $ 
 
  \noindent while $\cZ^j(0;\vs) = 0$ ; also, $0\leq \cZ^j(\tau)\leq 1$.

One shows that  it has a unique solution $\zz_s(\tau)=\zz_s^0(\tau)+\eps^2\zz_s^1(\tau,\eps)$  for small $\eps,$  where $\zz_s^0,\,\zz_s^1\sim 1$ and $\zz_s^0$  solves  the linear equation \eqref{kin eq}$|_{\eps=0 }$. Applying then the methods of  \cite{seclim} to the  estimates of Section 3, one would get

\smallskip\noindent{\bf Main theorem}.
{\it  Let $d\ge3,\,M\ge2$ be fixed. Then the energy spectrum  $\no_{s,L}(\tau)$  of the (M-order) quasisolution $\cA_s(\tau)$ of \eqref{ku5}, \eqref{in_cond} satisfies  $\no_{s,L}(\tau)\leq C_s$ and is $\eps^{2M}$-close to the solution $\zz_s (\tau)$ of mWKE \eqref{kin eq}. }
	
\medskip  Note that the result is deliberately formulated somewhat vaguely, since an elaboration of its details, along the lines of \cite{seclim} would need some extra  efforts (but no new ideas).\medskip 

We should also stress that this Main Theorem supposes the validity of our MC, moreover, explicit formulas for the energy spectrum of the $M$-order quasisolution (see \cite{seclim}) show that it corresponds to the case $N=M+1$ of MC, and respectively, is based on the $N=M+1$ case of Theorem 3.2.

 \end{document}